\documentclass[12pt]{article}
\usepackage{amsmath,amssymb,latexsym,epsfig,enumerate}
\usepackage{psfrag}
\usepackage[usenames]{color}\usepackage{graphicx}

\definecolor{brown}{cmyk}{0, 0.72, 1, 0.45}
\definecolor{grey}{gray}{0.5}

\newcommand{\ignore}[1]{}

\newcommand{\set}[1]{\left\{#1\right\}}

\def\cE{{\cal E}}

\def\E{\mathbb{E}}

\def\a{\alpha}

\def\e{\epsilon}

\def\G{\Gamma}

\def\th{\theta}

\def\n{\nu}

\def\om{\omega}

\def\Pr{\mathbb{P}}

\newtheorem{lemma}{Lemma}
\newtheorem{theorem}{Theorem}

\newcommand{\brac}[1]{\left( #1\right)}
\newcommand{\bfrac}[2]{\brac{\frac{#1}{#2}}}

\newcommand{\proofstart}{{\bf Proof\hspace{2em}}}

\newcommand{\proofend}{\hspace*{\fill}\mbox{$\Box$}}

\def\gnp{\mathbb{G}_{n,p}}
\newcommand{\beq}[1]{\begin{equation}\label{#1}}
\newcommand{\eeq}{\end{equation}}

\begin{document}
\title{{\bf On the non-planarity of a random subgraph}}
\author{Alan Frieze\thanks{
Department of Mathematical Sciences, Carnegie Mellon University,
Pittsburgh PA15213, U.S.A. Supported in part by NSF grant CCF1013110.}
\and Michael Krivelevich\thanks{ School of Mathematical Sciences,
Raymond and Beverly Sackler Faculty of Exact Sciences, Tel Aviv
University, Tel Aviv 69978, Israel. E-mail: krivelev@post.tau.ac.il.
Research supported in part by a USA-Israel BSF Grant and by a grant
from Israel Science Foundation.} } \maketitle

\begin{abstract}
Let $G$ be a finite graph with minimum degree $r$. Form a random
subgraph $G_p$ of $G$ by taking each edge of $G$ into $G_p$
independently and with probability $p$. We prove that for any
constant $\epsilon>0$, if $p=\frac{1+\epsilon}{r}$, then $G_p$ is
non-planar with probability approaching 1 as $r$ grows. This
generalizes classical results on planarity of binomial random
graphs.
\end{abstract}

\parindent 0in
{\bf AMS Classification:} 05C80, 05C10.
\section{Introduction}

Planarity is a fairly classical subject in the theory of random
graphs. Already Erd\H{o}s and R\'enyi in their groundbreaking paper
\cite{ER} stated (re-casting their statement in the language of
binomial random graphs) that a random graph $\gnp$ has a sharp
threshold for non-planarity at $p=1/n$ in the following sense: if
$p=c/n$ and $c<1$ then the random graph $\gnp$ is with high
probability (whp) planar, while for $c>1$ $\gnp$ is whp non-planar.
The Erd\H{o}s-R\'enyi argument for non-planarity had a certain
inaccuracy, as was pointed by \L uczak and Wierman \cite{LW}, who
explained how the probable non-planarity result can be obtained by
other means.

The aim of this paper is to generalize this classical non-planarity
result to a much wider class of probability spaces. All graphs
considered in this paper are finite. For a graph $G=(V,E)$ and
$0\leq p\leq 1$ we can define the random graph $G_p=(V,E_p)$ where
each $e\in E$ is independently included in $E_p$ with probability
$p$. When $G=K_n$, the complete graph on $n$ vertices, $G_p$ becomes
the binomial random graph $\gnp$.

Here is the main result of the present paper.

\begin{theorem}\label{th1}
Let $G$ be a finite graph with minimum degree $r$ and let
$p=\frac{1+\e}{r}$, where $\e>0$ is an arbitrary constant. Then
$$\Pr(G_p\text{ is planar})\leq \th_r$$
where $\lim_{r\to \infty}\th_r=0$.
\end{theorem}

\section{Proof of Theorem \ref{th1}}
Our proof rests in large part on the following simple consequence of Euler's formula.
\begin{lemma}\label{lem1}
Let $G=(V,E)$ be a planar graph with $n$ vertices and $m$ edges and girth $g$. Then
$$m\leq \frac{g(n-2)}{g-2}<n+\frac{2}{g-2}n.$$
\end{lemma}
\proofstart

Let $f$ be the number of faces of a planar embedding of $G$. Then we have
$$m=n+f-2\text{ and }2m\geq gf.$$
\proofend

\medskip

\noindent{\bf Remark.} In fact, a much stronger statement (in terms
of its consequences) is true: for any $\e>0$ and any integer $t$,
there exists $g=g(\e,t)$ such that any graph $G=(V,E)$ of average
degree at least $2+\e$ and of girth at least $g$ contains a minor of
the complete graph $K_t$. This was observed in particular by K\"uhn
and Osthus in \cite{KO03}. Indeed, by deleting repeatedly vertices
of degree 0 or 1 and paths of vertices of degree 2 in $G$ of length
at least $2/\e$ we keep average degree at least $2+\e$, and
eventually arrive at a subgraph $G'$ of $G$ of minimum degree at
least 2, in which every path of degree 2 vertices has length at most
$2/\e$. Contracting now these degree 2 paths produces a graph $G^*$
with still high girth (the girth went down by a factor at most
$2/\e$), but the minimum degree of $G^*$ is already at least 3. Then
applying the main result of \cite{KO03} to $G^*$  gives a large
complete minor in $G^*$, which corresponds to a large complete minor
in $G$. Alternatively, this can be derived directly from a result of
Mader \cite{Mad01}.

\medskip

Before proving Theorem \ref{th1} it will be instructive to prove it
for the special case where $G=K_n$ i.e. to show that if
$p=\frac{c}{n}$ where $c>1$ is a constant then $\gnp$ is non-planar
whp.

\subsection{$\gnp$}
The non-planarity of $\gnp$ is already known even for $c=1+\om
n^{-1/3}$ provided $\om\to\infty$ with $n$, see {\L}uczak, Pittel
and Wierman \cite{LPW}, see also \cite{NRR} for very accurate
results on the probability of planarity in the critical window
$p=(1+O(n^{-1/3}))/n$. The analysis for $c=1+o(1)$ is quite
challenging, but for constant $c>1$ it follows simply from some well
known facts. Let $G_1$ be the largest connected component of $\gnp$ (well
known to be whp the unique component of linear size for $c>1$, the
so called giant component). It is known, see e.g. \cite{B}, that whp
$$|V(G_1)|\sim xn\text{ and }|E(G_1)|\sim cn(2x-x^2)/2$$
where $x$ is the unique solution in $(0,1)$ to $x=1-e^{-cx}$. This
gives
$$c=1+\frac{x}{2}+\frac{x^2}{3}+\cdots$$
and so if $c=1+\e$, $\e>0$ and small, then $x=2\e-\frac83\e^2+O(\e^3)$.

Thus in this case, whp,
$$
\frac{|E(G_1)|}{|V(G_1)|}\sim \frac{c(2-x)}{2}=1+\frac{\e^2}{3}+O(\e^3).
$$
Next let $g_0=10/\e^2$. Then if $X$ denotes the number of cycles in $\gnp$ of
length at most $g_0$,
$$\E(X)\leq \sum_{k=3}^{g_0}n^kp^k\leq g_0c^{g_0}.$$
So, whp, there are fewer than $\ln n$ cycles of length at most
$g_0$. So, by removing at most $\ln n$ edges from $\gnp$ we obtain a
sub-graph $G_1'$ with girth higher than $g_0$. Now
$$\frac{|E(G_1')|}{|V(G_1')|}\sim \frac{|E(G_1)|}{|V(G_1)|}\sim 1+\frac{\e^2}{3}+O(\e^3)
>1+\frac{2}{g_0-2}$$
for small enough $\e$. Lemma \ref{lem1} implies that $G_1'$ and
hence $G_1$ are both non-planar. In fact, choosing a larger value of
$g_0$ and then recalling the remark following the proof of Lemma
\ref{lem1} shows that  $\gnp$ has with high probability an
arbitrarily large complete minor.

\subsection{Proof of Theorem \ref{th1}}
All asymptotic quantities are to be interpreted for $r\to\infty$ i.e. if we say
$\xi=\xi(r)=o(1)$ then we mean
that $\limsup_{r\to \infty}|\xi|=0$. This includes the notion of high probability.
I.e. if an event $\cE$ occurs with probability $1-\xi(r)$ where $\limsup_{r\to \infty}|\xi|=0$
then we say that $\cE$ ocurs whp.

{\bf Notation:} If $X$ is a set of edges and $A,B$ are disjoint sets of vertices, then
$E_X(A,B)$ is the
set of edges in $X$ with one endpoint in $A$ and one endpoint in $B$. Furthermore,
$E_X(A)$ is the
set of edges in $X$ with both endpoints in $A$. We let $e_X(A,B)=|E_X(A,B)|$ and
$e_X(A)=|E_X(A)|$.

Our strategy for proving Theorem \ref{th1} will be to prove the existence, whp,
of a sub-graph which has large
girth and sufficient edge density to apply Lemma \ref{lem1}. For this we will
need the following lemma:
\begin{lemma}\label{lem2}
Let $0<c_1,c_2<1$ be constants. Let $T=(V,E)$ be a tree on $n$
vertices with maximum degree $\Delta=r^{o(1)}$. Let $F\subseteq
\binom{V}{2}$ with $|F|= c_1nr$. Form a random subset $F_p$ of $F$
by choosing every edge of $F$ to belong to $F_p$ independently and
with probability $p=\frac{c_2}{r}$. Then the graph $G=T\cup F_p$ is
non-planar, with probability $1-O(r^{-1+o(1)})$.
\end{lemma}
\proofstart
Set
$$\alpha =\frac{c_1c_2}{72}\text{ and } A = \frac{10^{10}}{c_1^5c_2^6}\,.$$
Let
$$
V_0=\{v:d_F(v)\ge Ar\}\,,
$$
where $d_X(v)$ is the degree of vertex $v$ in the graph induced by $X\subseteq F$.

Clearly
$$|V_0|\le \frac{2c_1n}{A}.$$
Let $F'$ be the set of
edges from $F$ with at least one endpoint in $V_0$.

{\bf Case 1:} $|F'|\ge c_1nr/2$.\\
If
$e_F(V_0)\geq c_1nr/4$, then the Chernoff bound for the binomial distribution
implies that
with probability $1-e^{-\Omega(n)}$
we have
$$e_{F_p}(V_0)\geq\frac{c_1nr}{4}\cdot\frac{c_2}{r}(1-o(1))>\frac{c_1c_2n}{5}\geq
\frac{c_2A}{10}|V_0|>4|V_0|.$$
In which case, the subgraph induced by $E_{F_p}(V_0)$ forms a non-planar graph.
Hence we can
assume from now on that $F'$ has at least $c_1nr/4$ edges with at most one endpoint
in $V_0$.

Define
$$
U_0 =\set{v\not\in V_0: d_{F'}(v)\ge \frac{c_1r}{8}}\,.
$$
Then
$$e_F(U_0,V_0)\geq |F'|-e_F(V_0)-\frac{c_1nr}{8}\geq \frac{c_1nr}{8}.$$
Now the Chernoff bound implies that with probability
$1-e^{-\Omega(n)}$ we have
$$e_{F_p}(U_0,V_0)\geq\frac{c_1nr}{8}\cdot\frac{c_2}{r}(1-o(1))\ge
\frac{c_1c_2n}{9}.$$
So if $|U_0|\leq \a n$ then
$$\frac{e_{F_p}(U_0,V_0)}{|U_0|+|V_0|}\geq \frac{c_1c_2n/9}{c_1c_2n/72+2c_1n/A}\geq 4.$$
In which case, the subgraph induced by $E_{F_p}(V_0,U_0)$ forms a
non-planar graph.

If $|U_0|\ge \alpha n$, then define a (random) subset $W_0$ by:
$$
W_0=\{v\in U_0: d_{F_p}(v,V_0)\ge 5\}\,.
$$
The distribution of $|W_0|$ dominates $Bin(|U_0|,q_1)$ where
$$q_1=\Pr(Bin(c_1r/8,c_2/r)\geq 5)\geq \binom{\frac{c_1r}{8}}{5}
\left(\frac{c_2}{r}\right)^5 \brac{1-\frac{c_2}{r}}^{c_1r/8}\geq
(1-o(1))\left(\frac{c_1c_2}{8}\right)^5
\cdot\frac{e^{-c_1c_2/8}}{120}.$$ The Chernoff bounds then imply
that with probability $1-e^{-\Omega(n)}$ we have
$$|W_0|\ge \frac{|U_0|q_1}{2} \ge \frac{\a nq_1}{2}\,,
$$
and by definition $e_{F_p}(W_0,V_0)\ge 5|W_0|$. This is at least $4|W_0\cup V_0|$.
In which case, the subgraph induced by $E_{F_p}(V_0,W_0)$ forms a non-planar graph.
This completes the analysis for Case 1.

{\bf Case 2:} $|F'|\le c_1nr/2$.\\
Define
$F''=F\setminus F'$ and observe that $|F''|\ge c_1nr/2$, and by definition the
maximum degree of $F''$ is at most $Ar$.

Observe that $T\cup F_p''$ has
with probability $1-e^{-\Omega(n)}$
at least
$$n-1+\frac{c_1nr}{2}\cdot\frac{c_2}{r}(1-o(1))> n(1+\epsilon)$$
edges for some positive $\epsilon=\epsilon(c_1,c_2)$. It thus
suffices to show that the number of ``short`` cycles in $T\cup
F''_p$ is $o(n)$ whp, and then to use Lemma~\ref{lem1}.

For constants $\ell,t=O(1)$ let us estimate the expected number of
cycles of length $\ell$ in $T\cup F_p''$ having $t$ edges from $T$.
We choose an initial vertex $v$ in $n$ ways, then decide about the
placement of the edges of $T$ in the cycle in $O(1)$ ways. We thus
get a sequence $P_1*P_2\cdots * P_{t+1}$, where the stars correspond
to edges from $T$ and $P_i$ is a path of length $\ell_i$ in $F''$
for $i=1,2,\ldots,t+1$. Now, a path of length $\ell_i$ in $F''$,
starting from a given point, can be chosen in at most
$(Ar)^{\ell_i}$ ways, an edge from $T$ from a given vertex can be
chosen in at most $\Delta(T)=r^{o(1)}$ ways, and finally the last
path of length $\ell_{t+1}$, connecting two already chosen vertices,
can be chosen in at most $(Ar)^{\ell_{t+1}-1}$ ways. Altogether the
number of such cycles in $T\cup F''$ is $n\cdot r^{o(1)}\cdot
O\left(r^{\ell_1+\ldots+\ell_{t+1}-1}\right)$. The probability
for such a cycle to survive in $T\cup F_p''$ is
$O\left(r^{-(\ell_1+\ldots+\ell_{t+1})}\right)$. We thus expect
$O(n/r^{1-o(1)})$ such cycles. Summing over all choices of $\ell$
and $t$ we get that the expected number of cycles of length $O(1)$
in $T\cup F_p''$ is $O(n/r^{1-o(1)})$, and thus the Markov
inequality implies we get fewer than $n/\ln r$ cycles, with
probability $1-O(r^{-1+o(1)})$. By choosing $\ell$ sufficently large
and deleting one edge from each cycle of length at most $\ell$, we
get a graph of large constant girth, with $n$ vertices and at least
$(1+\epsilon/2)n$ edges --- which is non-planar, by Lemma
\ref{lem1}. \proofend

We now set about using the above lemma.
We let $G_p=G_1\cup G_2$ where $G_i=G_{p_i},i=1,2$, and $p_1=\frac{1+\e/2}{r}$ and
$(1-p_1)(1-p_2)=1-p$ so that
$p_2=\frac{\e+O(\e^2)}{2r}$.

\subsection{Proof outline}
Before going to concrete details we provide a short outline of the
proof. We start by probing relatively few vertices and their
incident random edges till we find a vertex $v$ whose degree in
$G_1$ is at least $d=\ln^{1/2}r$.The immediate neighborhood of $v$
in $G_1$ is large enough to support the growth of (some version of)
the BFS tree $T_k$ from $v$ until it accumulates about $i_0=\ln^3r$
vertices, while its frontier $S_k$ is of size $\Theta(\e)|T_k|$.

{From} this point on, we proceed iteratively, at each iteration
looking at the current tree $T_k$, its frontier $S_k$ and the edges
of $G$ touching $S_k$. If many of these edges go back to $T_k$, we
can sprinkle them in $G_2$ and apply Lemma \ref{lem2} to argue that
the resulting random graph is whp non-planar. Otherwise, many of the
edges touching $S_k$ leave $T_k$, which allows us to expose them in
$G_1$ and to add yet another layer of substantial size to the
current tree, while controlling its maximum degree, and to proceed
to the next iteration. This growth process cannot go forever, as $G$
is finite, and thus it eventually collapses, with the first
alternative above being applicable, thus resulting in a non-planar graph
whp.

\subsection{Initial Tree Growth}\label{itg}
We begin by repeatedly choosing a vertex $v\in V$ and analysing a
restricted breadth search (RBFS) from $v$ until we succeed in obtaining a certain
condition, see \eqref{eq4} below. Basically, we need to find $v$ which has
sufficiently many neighbors in $G_1$.
So, let $S_0=\set{v}$.
In general let
$$S_{i+1}=\bigcup_{w\in S_i}(RN(w)\setminus (T_i\cup B))$$
where
\begin{itemize}
\item $B, |B|=o(r)$, is a set of vertices that already been rejected by our search.
\item $T_i=\bigcup_{j=0}^iS_j$.
\item $RN(w)$ denotes the first $Bin(r_1,p_1)$, $r_1=r-O(i_0\ln r)$, neighbors of $w$
in $G_1$, where
$$i_0=\ln^3r,$$
By first we assume that $V(G)=[n]$ for some integer $n$.
Then we mean that we try the first $r_1$ $G$-neighbors of a vertex $w$ in numerical
value to see if they are neighbors
of $w$ in $G_1$. The edges found will be part of a subgraph $H_1$ and we only keep
the first edge found to each vertex
added. In this way, $H_1$ will be a tree.
\end{itemize}
Our initial aim in RBFS is to find a smallest $k$ such that
\beq{eq4} i_0\leq |T_k|\leq 2i_0\text{ and }\frac{|S_k|}{|T_k|}\in
\left[\frac{\e}{4}, \frac{3\e}{4}\right]. \eeq Let $d=\ln^{1/2}r$.
We first look for $v$ such that $d\leq |S_1|\leq \ln^2r$. This is
quite simple. Let $l_0=(2d)^d=o(r)$ and suppose that we have already
examined $v_1,v_2,\ldots,v_l,\,l\leq l_0$, without success. We
choose $v\notin B_l=\set{v_1,v_2,\ldots,v_l}$ and examine the first
$r-o(r)$ neighbours of $v$ that are not in $B_l$. The probability
that $v$ has at least $d$ neighbors in $G_1$ is greater than
$\binom{r-o(r)}{d}p_1^d(1-p_1)^{r-o(r)}\geq d^{-d}$. So, the
probability we have not found $v$ with large enough degree after
$l_0$ trials is less than $(1-d^{-d})^{l_0}=o(1)$. Furthermore, the
probability $v$ has more than $\log^2r$ neighbors is less than
$\binom{r}{\ln^2r}p_1^{\ln^2r}\leq e^{-\ln^2r}$. We can therefore
assume that we can find a suitable $v$ with $d\leq |S_1|\leq
\ln^2r$, where $B=B_l$ is of size $o(r)$.

Suppose now that $S_i,T_i,i\geq 1$ do not satisfy \eqref{eq4} and that $|T_i|\leq 2i_0$.
We observe first that the distribution of the size of $S_{i+1}$ is dominated by
$Bin((r-o(r))|S_i|,p_1)$. In fact we bound $|S_{i+1}|$ from above by the number of
edges from $S_i$ to $S_{i+1}$.
We examine the first $r_1=r-o(r)$ $G$-neighbors of each $v$ in $S_i$
and include an edge $vw$ in our count if the edge $vw$ is in $G_1$.
Therefore
\beq{B1}
\Pr(|S_{i+1}|\geq(1+2\e/3)s\mid |S_i|=s)\leq e^{-\Omega(\e^2s)}.
\eeq
We can also argue that $|S_{i+1}|$ dominates a binomial $Bin(|S_i|(r-o(r)),p_1)$.
The $o(r)$ term here differs from the one used in the upper bound.
We will have to exclude edges to those $G$-neighbors that have already been placed
in $S_{i+1}$ and to those
$G$-neighbors in $B_l$. Because we are looking for a lower bound which
is less than $i_0$, we can claim to get at least the result of $|S_i|(r-o(r))$ trials
with success probability $p_1$. Therefore
\beq{B2}
\Pr(|S_{i+1}|\leq(1+\e/3)s\mid |S_i|=s)\leq e^{-\Omega(\e^2s)}.
\eeq
So we can assume that $\a_1\geq d$ and $|S_i|/|S_{i-1}|=\a_i\in [(1+\e/3),(1+2\e/3)]$ for
$i\geq 2$. And then
\beq{ST}
\frac{|S_i|}{|T_i|}=\frac{\a_1\a_2\cdots\a_i}{1+\a_1+\a_1\a_2+\cdots+\a_1\a_2\cdots\a_i}.
\eeq
The expression \eqref{ST} is minimised (resp. maximised) by putting $\a_i=(1+\e/3)$
(resp. $=(1+2\e/3)$) for
$i\geq 2$.
It follows that whp
\beq{eq99}
\frac{|S_i|}{|T_i|}= \frac{\a_1\th^{i-1}(\th-1)}{\th-1+\a_1(\th^i-1)}
\eeq
for some $\th\in [(1+\e/3),(1+2\e/3)]$.

Thus we will achieve \eqref{eq4} whp. Here we use two facts: (i) the sum of the failure
probabilities in \eqref{B1},\eqref{B2} is bounded by
$\sum_{s\geq d}e^{-\Omega(\e^2s)}=o(1)$; (ii) We have assumed that $|S_1|=o(i_0)$ which
means that the value of $k$ in \eqref{eq4} is $\omega(1)$ which in turn means that we
need only consider large $i$ in
\eqref{eq99}. In which case the ratio in \eqref{eq99} is asymptotically
$\frac{\th^{-1}(\th-1)}{\th^i}$.

\subsection{Remaining Tree Growth}
Let us consider the current tree $T_k$, which is of size $\Omega(i_0)$, and its frontier
$S_k$ of
size $|S_k|=s_k=\Theta(\epsilon |T_k|)$. Choose $r_1=r-o(r)$
arbitrary edges incident to each vertex of $S_k$, denote the
obtained set by $E_k$, $|E_k|\le rs_k$. If $E_k$ has $\Theta(rs_k)$
edges inside $V(T_k)$, then sprinkling the edges of $E_k$  with
probability $p_2$ produces whp a non-planar graph on $V(T_k)$ by Lemma \ref{lem2}.

We can therefore assume that $E_k$ has at least
$(1-\frac{\epsilon}{10})rs_k$ edges between $S_k$ and $V\setminus T_k$.

Let $V_0=\{v\not\in T_k: d_{E_k}(v,S_k)\ge r\ln r\}$. Clearly,
$|V_0|\le s_k/\ln r$. If $E_k$ has at least
$\epsilon rs_k/10$ edges between $S_k$ and $V_0$, then in the
random subset of $E_k$, formed by taking each edge independently and
with probability $p_1$, there is whp a set $W_0$ of
$|W_0|=\Theta(s_k)$ vertices $v\in S_k$, whose degrees $\eta_v$ into $V_0$ are
at least three. Indeed, there will be at least
$\e s_k/20$ vertices $S_k'$ in $S_k$ that have at least $\e r/20$ neighbours in
$V_0$.
Each vertex in $S_k'$ has a probability of at least $\xi=
\binom{\e r/20}{3}p_1^3(1-p_1)^{\e r/20-3}\ge \e^3 10^{-10}$
of having $\eta_v\geq 3$, and these events are independent.
Thus whp $|W_0|\geq |S_k'|\xi/2$ and the bipartite subgraph
of $G_p$ induced by $W_0,V_0$ has more than $2(|W_0|+|V_0|)$ edges and so is non-planar.
We can assume therefore that $E_k$ has at least
$(1-\frac{\epsilon}{5})rs_k$ edges between $S_k$ and $V_1=V\setminus(T_k\cup V_0)$.
Denote this set of edges by $F_k$.

Form a random subgraph $R_k$ of $F_k$ by taking each edge
independently and with probability $p_1$.
\begin{enumerate}[{\bf P1}]
\item Then the Chernoff bound implies that with probability
$1-e^{-\Omega(s_k)}$, $|R_k|\ge
(1+\frac{\epsilon}{5})s_k$.
\item Furthermore, we will show next that with probability $1-\e_1(r)$ at most
$2s_k/\ln r$ of these edges are incident with vertices in
$V_2\subseteq S_k$ whose degree in $R_k$ is more than $\ln\ln r$.
\end{enumerate}
The value of $\e_1(r)=\e_1'(r)+\e_1''(r)$ is obtained from \eqref{e11} and
\eqref{e1} below. Indeed, if $v\in S_k$ then
$$\Pr(d_{R_k}(v)\geq \ln\ln r)\leq \Pr(Bin(r,p_1)\geq \ln\ln r)\leq
q_2=\bfrac{2e}{\ln\ln r}^{\ln\ln r}$$ and
$$\Pr(d_{R_k}(v)\geq \ln r)\leq \Pr(Bin(r,p_1)\geq \ln r)\leq
q_3=\bfrac{2e}{\ln r}^{\ln r}.$$
Thus the number of edges in $R_k$ that are incident with $v\in V_2$ is bounded by\\
$Bin(s_k,q_2)\ln r+Bin(s_k,q_3)r$. We observe that because $s_k\geq
\e i_0/8\gg \ln^2r$ we can write \beq{e11} \Pr(Bin(s_k,q_2)\geq
s_k/\ln^2r)\leq \e_1'(r)=(e\ln^2r\,q_2)^{s_k/\ln^2r} \eeq and
\beq{e1} \Pr(Bin(s_k,q_3)\geq s_k/r^2)\leq \e_1''(r)=
\begin{cases}r^2q_3&0\leq s_k<r^3\\ (er^2\,q_3)^{s_k/r^2}&s_k\geq r^3\end{cases}.
\eeq

Let $N_k$ be the set of neighbors of $S_k$ defined by edges in $R_k$.
We observe that $|N_k|$ is the sum of independent
Bernouilli random variables. We consider two cases depending
on the value of $\E(|N_k|)$ w.r.t. the random set $R_k$.
Splitting the argument this way will not condition $R_k$
or $N_k$.

{\bf Case 1:} $\E(|N_k|)\ge\brac{1+\frac{\epsilon}{10}}s_k$.

We first observe that
\beq{M1}
\Pr\brac{|N_k|\leq \brac{1+\frac{\epsilon}{20}}s_k}\le
e^{-\e^2s_k/1000}.
\eeq
We therefore assume that
$$|N_k|\geq \brac{1+\frac{\epsilon}{20}}s_k.$$
$R_k$ contains a subset $R_k'$ of size
$\n_k=\brac{1+\frac{\epsilon}{25}}s_k$ such that the degrees of all
the vertices in $S_k$ w.r.t. $R_k'$ are at most $\ln\ln r$, and
every vertex outside $T_k$ is incident to at most one edge from
$R_k'$ and there are $\n_k$ vertices outside $T_k$ incident to an
edge in $R_k'$. We obtain this by removing edges incident with $V_2$
and by then deleting edges incident with $N_k$ to get degree at most
one. Use $R_k'$ to form the next frontier of size
$(1+\frac{\epsilon}{25})s_k$, composed of the endpoints of the edges
of $R_k'$ outside $S_k$; proceed to the next round.

{\bf Case 2:} $\E(|N_k|)\le\brac{1+\frac{\epsilon}{10}}s_k$.

\begin{description}
\item[{\bf Q1}] $|N_k|\le (1+\frac{\epsilon}{8})s_k$ whp. Indeed,
$$\Pr\brac{|N_k|\geq \brac{1+\frac{\epsilon}{8}}s_k}\le
e^{-\e^2s_k/5000}.$$
\item[{\bf Q2}$\equiv${\bf P1}]
$|R_k|\ge (1+\frac{\epsilon}{5})s_k$ whp. Indeed, $|R_k|=Bin(|F_k|,p_1)$ and so
$$\Pr\brac{|R_k|\le \brac{1+\frac{\epsilon}{5}}s_k}\leq e^{-\e^2s_k/1200}$$
for small $\e>0$.
\item [{\bf Q3}]There are $o(s_k)$ short cycles in $T_k\cup R_k$ whp.
For this calculation we consider the graph $\G_k$ induced
by the edges in $E(T_k)\cup F_k$. This has vertex set
$V(T_k)\cup N_k$.
Here the
expectation calculation is quite similar to that of the lemma. We
use the fact that $V_0$ has been excluded, and therefore all
relevant vertices outside of $T_k$ have their degrees into $S_k$
bounded by $r\ln r$. Also all degrees in $T_k$ are $r^{o(1)}$ by our
construction.\\
{\bf Details:} For constants $\ell,t=O(1)$ let us estimate the
expected number of cycles of length $\ell$ in $\G_k$ having $t$
edges from $T_k$. We choose an initial vertex $v$ in $O(s_k)$ ways,
then decide about the placement of the edges of $T_k$ in the cycle
in $O(1)$ ways. We thus get a sequence $P_1*P_2\cdots * P_{t+1}$,
where the stars correspond to edges from $T_k$ and $P_i$ is a path
of length $\ell_i$ for $i=1,2,\ldots,t+1$ using edges in $F_k$. Now,
a path of length $\ell_i$ using edges in $F_k$, starting from a
given point, can be chosen in at most $(r\ln r)^{\ell_i}$ ways, an
edge from $T_k$ from a given vertex can be chosen in at most
$\Delta(T_k)=r^{o(1)}$ ways, and finally the last path of length
$\ell_{t+1}$, connecting two already chosen vertices, can be chosen
in at most $(r\log r)^{\ell_{t+1}-1}$ ways. Altogether the number of
such cycles in $\G$ is $s_k\cdot r^{o(1)}\cdot
\tilde{O}\left(r^{\ell_1+\ldots+\ell_{t+1}-1}\right)$ ways. The
probability for such a cycle to survive in $T_k\cup R_k$ is
$O\left(r^{-(\ell_1+\ldots+\ell_{t+1})}\right)$. We thus expect
$O(s_k/r^{1-o(1)})$ such cycles. Summing over all choices of $\ell$
and $t$ we get that the expected number of cycles of length $O(1)$
in $T\cup F_p''$ is $O(s_k/r^{1-o(1)})$, and thus the Markov
inequality implies we get fewer than $s_k/\ln r$ cycles, with
probability $1-O(r^{-1+o(1)})$.
\end{description}
By choosing $\ell$ sufficiently large and removing edges from the short cycles
(length $\le\ell$) leaves a graph of average degree
$2+\Theta(\epsilon)$ and without short cycles. This is non-planar
by Lemma \ref{lem1}.

As a final note in proof, we argue about the probability that this construction fails.
We have seen that the initial tree growth in Section \ref{itg} succeeds whp.
The success of
the remaining tree growth rests on the probabilities in {\bf P1,P2} being high enough.
These events
need to happen multiple times, whereas other events are only required to occur once.

For {\bf P1} and \eqref{M1} we verify that $\sum_{t\geq i_0}e^{-\Omega(t)}=o(1)$ and for
{\bf P2} we verify that
$$\sum_{t\geq i_0}(e\log^2r\,q_2)^{-t/\log^2r}+\sum_{t=0}^{r^3}r^2q_3+
\sum_{t\geq r^3}(er^2\,q_3)^{t/r^2}=o(1).$$

\section{Concluding remarks}
We have proven that for every finite graph $G$ of minimum degree
$r\gg1$, a random subgraph $G_p$ of $G$, with
$p=p(r)=\frac{1+\epsilon}{r}$ and $\epsilon>0$ being an arbitrary small
constant, is whp non-planar. This generalizes the classical
non-planarity results for binomial random graphs $\gnp$. It should
be noted that for a statement of such generality we cannot hope to
have a matching lower bound on $p(r)$. Indeed, if $G$ is a
collection of, say, $2^{r^3}$ vertex disjoint cliques $K_{r+1}$,
then for any constant $c>0$, the random subgraph $G_p$, $p=c/r$,
retains whp one of the cliques $K_{r+1}$ in full and is thus whp
non-planar.

Notice that our proof, with fairly straightforward and simple
adjustments,  shows in fact that under the conditions of Theorem
\ref{th1} the random subgraph $G_p$ is typically not only
non-planar, but has a complete minor of arbitrarily large constant
size. This can be obtained by employing the remark following Lemma
\ref{lem1}. It would be interesting to determine the largest
$t=t(r)$ such that under the same conditions the random graph $G_p$
has whp a minor of a complete graph $K_t$. For the case of binomial
random graphs $\gnp$  Fountoulakis, K\"uhn and Osthus showed
\cite{FKO1} that for any $c>1$, the random graph $\gnp$ with $p=c/n$
has whp a complete minor of order $\sqrt{n}$. (See also \cite{KS}
for results for other values of $p=p(n)$, and \cite{FKO2} for
results on random regular graphs and for $\gnp$ in the slightly
supercritical regime).

The main theorem of this paper can be viewed as yet another
contribution to a growing sequence of results about properties of
random subgraphs of graphs of given minimum degree. We can mention
here \cite{KS2}, who showed that if $G$ is a finite graph of minimum
degree $r$ and $p=\frac{1+\e}{r}$, then the random graph $G_p$
contains whp a path of length linear in $r$, and also \cite{KLS},
where it is proven that under the same assumptions on the base graph
$G$ and when taking $p=\frac{(1+o(1))\ln r}{r}$, the random graph
$G_p$ contains whp a path of length at least $r$, in both cases
substantially generalizing classical results about binomial random
graphs. One can certainly anticipate more results of this type to
appear in the near future.

{\bf Acknowledgement:} We thank the referees for a careful reading.

\end{document}